\documentclass[a4paper]{article}
\pdfoutput=1
\usepackage{graphicx}
\usepackage{amsmath,bbm,amssymb,array,verbatim}
\usepackage{hyperref}
\usepackage{algorithm}
\usepackage[noend]{algpseudocode}

\title{e-Valuate: A Two-player Game on Arithmetic Expressions}
\author{Sarang Aravamuthan\thanks{Ignite R\&D Labs, Tata Consultancy Services,  Chennai, India. \newline
\indent \hspace{1.5mm} \texttt{sarang.aravamuthan@tcs.com, b.ganguly@tcs.com}} \and Biswajit Ganguly\footnotemark[1]}

\newtheorem{thm}{Theorem}

\newcommand{\game}{\mbox{\emph{game}}}

\newcommand{\E}{\mbox{$\mathbb{E}$}}
\newcommand{\F}{\mbox{$\mathbb{F}$}}
\newcommand{\N}{\mbox{$\mathbb{N}$}}

\newcommand{\val}{\mbox{\rm val}}
\newcommand{\minimax}{\mbox{\emph{minimax}}}

\newcommand{\tree}{\mbox{\emph{tree}}}
\newcommand{\real}{\hbox{I\kern-.2em\hbox{R}}}
\newcommand{\qed}{\hfill \ensuremath{\Box}}

\algrenewcommand\algorithmicthen{}
\algrenewcommand\algorithmicdo{}
\algrenewcommand\algorithmiccomment[1]{\hfill\(\triangleright\) \texttt{#1}}%
\algdef{SE}[FUNCTION]{Function}{EndFunction}%
   [2]{\algorithmicfunction\ {#1}\ifthenelse{\equal{#2}{}}{}{(#2)}}%
   {\algorithmicend\ \algorithmicfunction}%

\begin{document}
\maketitle

\begin{abstract}
\noindent
e-Valuate is a game on arithmetic expressions. The players have contrasting roles of maximizing and minimizing the given expression. The maximizer proposes values and the minimizer substitutes them for variables of his choice. When the expression is fully instantiated, its value is compared with a certain minimax value that would result if the players played to their optimal strategies. The winner is declared based on this comparison.
\smallskip

We use a game tree to represent the state of the game and show how the minimax value can be computed efficiently using backward induction and alpha-beta pruning. The efficacy of alpha-beta pruning depends on the order in which the nodes are evaluated. Further improvements can be obtained by using transposition tables to prevent reevaluation of the same nodes. We propose a heuristic for node ordering.  We show how the use of the heuristic and transposition tables lead to improved performance by comparing the number of nodes pruned by each method.
\smallskip

We describe some domain-specific variants of this game. The first is a graph theoretic formulation wherein two players share a set of elements of a graph by coloring a related set with each player looking to maximize his share. The set being shared could be either the set of vertices, edges or faces (for a planar graph). An application of this is the sharing of regions enclosed by a planar graph where each player's aim is to maximize the area of his share. Another variant is a tiling game where the players alternately place dominoes on a $8 \times 8$ checkerboard to construct a maximal partial tiling. We show that the size of the tiling $x$ satisfies $22 \le x \le 32$ by proving that any maximal partial tiling requires at least $22$ dominoes.

\smallskip

\noindent
{\bf Keywords:} Arithmetic expressions, game tree, games on graphs, tiling 
\end{abstract}

\section{Introduction}
\label{sec:intro}

Given an arithmetic expression $\E$ involving variables and the standard operators ($+, -, *$ and $/$), players MAX and MIN evaluate $\E$ with contrasting goals; MAX would like to maximize $\E$ while MIN would like to minimize $\E$. Towards this end, they take turns to instantiate the variables. MAX starts and, at each move, proposes a value (digit 0--9) that MIN substitutes for a variable of his choice. When the expression is fully instantiated, it is evaluated and compared with a certain minimax value that would result if the players played to their optimal strategies. Let $\val (\E)$ be the value of $\E$ at the end of the game and $\minimax(\E)$ be the minimax value. The winner is then determined in the following way.
\begin{itemize}
\item If $\val(\E) > \minimax (\E)$ then MAX is declared the winner.
\item If $\val(\E) < \minimax (\E)$ then MIN is declared the winner.
\item If $\val (\E) = \minimax (\E)$, then the game is a draw.
\end{itemize}
\smallskip

\noindent
For example, if $\E = X*(Y - Z)$, a possible sequence of moves is
\begin{enumerate}
\item MAX chooses 5 and MIN replaces $X$ with 5 so that $\E = 5*(Y-Z)$.
\item Next MAX chooses 3 that MIN substitutes 3 for $Z$ leading to $\E = 5*(Y-3)$.
\item Finally MAX chooses 9 which MIN substitutes for the remaining variable $Y$ and the final value for the expression is $5(9-3) = 30$.
\end{enumerate}

\noindent
With more strategic play from either player, the expression is evaluated differently. 
For instance, with the same moves from MAX and optimal play from MIN, 
the substitutions would be $5\rightarrow Y, 3\rightarrow X$ and $9\rightarrow Z$ and the expression evaluates to $-12$. With optimal play from both players, a possible sequence of moves is $6\rightarrow Y, 3\rightarrow X$ and $0\rightarrow Z$ with $\E$ evaluating to the minimax value 18. 
\smallskip

We will refer to this version of the game as \emph{e-Valuate}. Specific instances of the game have appeared in books on mathematical puzzles. For example, in \cite{Wink}, the expression is a difference of two four digit numbers and the reader is asked to find the minimax value.
\smallskip

\noindent
Some possible variations on this form of the game are the following.
\begin{itemize}
\item The expression as well as the domain can be generalized. For example, other mathematical operators can be introduced in the expression and the domain can include other values over which the expression can be evaluated.
\item An alternate way of playing the game is for the players to switch roles at the end of the game and reevaluate the expression. If the expression evaluates to a larger value in one of the games then the maximizer in that game is the winner. This version could be applicable when the number of variables is large enough that computing the minimax value is infeasible.
\end{itemize}
Another variant is for the first player to take on the role of the minimizer and the second player that of the maximizer. This is however equivalent to the original version since $\min (\E) = -\max (-\E)$ and $\max(\E) = -\min(-\E)$ where the minimum and maximum are carried out over the domain of the variables. Thus, the final value under optimal play from both players is $-\minimax (-\E)$.
\smallskip

Minimax is a more general term and applies to any two player zero-sum game \cite{RN}. By using a \emph{game tree} to represent the states of the game and the moves of the players, the minimax algorithm can be used to determine the best move at each position in the game in the following manner. First values are assigned to the leaf nodes using an \emph{evaluation function}. Next, the players MAX and MIN attempt to maximize and minimize the value of the nodes corresponding to their turn of play. For an intermediate node that corresponds to MAX's turn to play, the value of the node is the maximum of the values of its children. Similarly, for an intermediate node that corresponds to MIN's turn to play, the value of the node is the minimum of the values of its children. The value at the root is the minimax value of the game. For example, if a game is designed such that under optimal play, MIN has a winning strategy, and the leaf nodes are assigned a value of $+1$ or $-1$ according to whether the corresponding position is a win for MAX or MIN, then the minimax value will be $-1$.
\smallskip

Several optimizations to this method of computing the minimax value have been studied \cite{IWJF, RN}. Some well known techniques are
\begin{itemize}
\item \emph{Alpha-beta pruning}: This is a windowing procedure that starts with an interval of $(-\infty, +\infty)$ for the minimax value. As nodes are evaluated, the window shrinks and any node that evaluates to a value outside this window is pruned along with the subtree rooted at that node.
\item \emph{Negascout}: Negascout \cite{Rei} works by assuming that for each node, the first child will be in the \emph{principal variation} (the sequence of moves leading to the minimax value). It uses a null search window for the remaining children and on failure, uses a full search window. Thus this method is most effective when there is a good ordering for evaluating the nodes.
\item \emph{Transposition tables}: This is a memoization technique where the values of nodes that are evaluated are stored and retrieved when another node that corresponds to the same game position has to be evaluated. This effectively prunes the subtree rooted at that node.
\end{itemize}
\smallskip

The computational challenge in e-Valuate is an efficient way of determining the minimax value in order to identify the winner. We show how these techniques lead to more efficient ways of determining $\minimax(\E)$. 
\smallskip

The game can also be cast in other domains. In particular we study some graph theoretic versions of the game wherein two players share a set of entities of a graph by coloring a different set with each player aiming to maximize his share. The entities being shared can be the set of vertices, edges or faces of the graph.
\smallskip

In the next section, we introduce the game tree for e-Valuate and show how the minimax value can be computed using \emph{backward induction}. In Section~\ref{sec:alpha}, we show how improved performance can be obtained by combining the minimax algorithm with alpha-beta pruning. We describe these methods in the context of our game. The efficacy of alpha-beta pruning methods depends on the order in which the the children of each node are evaluated. We describe a heuristic for determining this order. Further improvements can be obtained by avoiding repeated reevaluation of the same game position through the use of transposition tables. In Section~\ref{sec:Imp}, we provide implementation details and compare the number of nodes pruned by the two methods, alpha-beta and alpha-beta with node ordering and transposition tables, for different arithmetic expressions. In Section~\ref{sec:G} we introduce some versions of e-Valuate that can be played on graphs. In Section~\ref{sec:T} we describe another variant that is based on a maximal partial tiling of an $8 \times 8$ checkerboard by dominoes. In particular we show that any maximal partial tiling requires at least $22$ dominoes. We conclude with some unanswered questions related to this game.
\smallskip

We fix some notations. For an arithmetic expression $\E$, let $n$ be the number of variables in $\E$. $\E (i \rightarrow X)$ denotes the expression $\E$ with variable $X$ replaced by $i$. For a domain $D$ different from the set of digits, we will denote by $\minimax(\E, D)$, the minimax value of $\E$ where each variable in $\E$ takes values from $D$. The set $[k] = \{1, \dots, k\}$. All the games we describe in this paper are sequential two-player games played by MAX and MIN with MAX making the first move.
\smallskip

For general aspects of game theory, see \cite{OR}; \cite{net} is a useful online resource for lectures, glossary of terms and articles related to game theory. The game algorithms we have outlined above are well documented in books on artificial intelligence (e.g. \cite{IWJF}, \cite{RN}). \cite{BM} provides an excellent introduction to graph theory.

\section{The Game Tree for e-Valuate}	
\label{sec:Games}

In the framework of game theory, e-Valuate can be classified as a finite, sequential, two person game with perfect information. It is finite as the game ends after a finite number of moves, sequential since the players take turns in making their moves (rather than move simultaneously as, say, in the \emph{rocks, paper and scissors} game) and it's a game of perfect information as each player is aware of the other's moves at any point in the game.
\smallskip

Sequential games with perfect information can be represented using a game tree. The root of the tree corresponds to the initial configuration of the game (in our case, the expression $\E$) and the edges represent possible moves that the players make. Each node in the tree represents a position in the game. The root and the leaf nodes are MAX nodes and the nodes at intermediate levels are alternately MAX and MIN nodes and represent positions where the MAX or MIN has to make a move. Thus each MAX node has 10 children corresponding to 10 possible moves (choosing any digit). A MIN node at a height $d$ has $(d+1)/2$ children that correspond to $(d+1)/2$ uninstantiated variables. The height of the tree is $2n$. We will denote by $\tree(\E)$, the game tree corresponding to $\E$. 
\smallskip

The number of nodes in the tree, $T(n)$, depends only on $n$ and satisfies the recursion
\begin{equation}
T(n) = 11 + 10nT(n-1)
\label{eqn:1}
\end{equation}
which follows from observing that the root node has 10 children each of which has $n$ children that correspond to game trees on expressions with $(n-1)$ variables.
\smallskip

We can use this to bound $T(n)$ by
\[ 2n! 10^n \le T(n) \le 2n! 10^n e^{1/10} \]
from the following argument. Let $N = n! 10^n$ be the number of leaves of $\tree(\E)$. Starting from the bottom and counting the number of nodes at each level we get
\begin{eqnarray*}
T(n) & = & N + N/1 + N/(1*10) + N/(1*10*2) + \dots + N/(n! 10^n) \\
&= & \sum_{i=0}^n N/(i! 10^i) + \sum_{i=1}^n N/(i! 10^{i-1}) \le 2  \sum_{i=0}^n N/(i! 10^i) \le 2e^{1/10} N
\end{eqnarray*}
as desired.
\smallskip

\noindent
We identify each node in a tree by
\begin{itemize}
\item a sequence of instantiations of the variables and possibly an additional digit (for a MIN node). For example if $\E = (10 - X)*Y$, then a MAX node in $\tree (\E)$ is $\{ 1 \rightarrow Y\}$ and a MIN node is $\{ 1 \rightarrow Y, 3\}$.  Thus MAX nodes correspond to partially instantiated expressions and MIN nodes to (expression, digit) pairs.
\item a value which is the minimax value of the partially instantiated expression for a MAX node and the minimum of the minimax values of the children for a MIN node. This is the value $\E$ would evaluate to under optimal play starting from the position given by the node. This is also referred to as the \emph{score} of the position given by the node \cite{IWJF}.
\end{itemize}

\smallskip

\begin{figure}[htbp]
\begin{center}
\includegraphics[height=80mm, width=1.2\textwidth,trim = 25mm 5mm 0mm 60mm, clip]{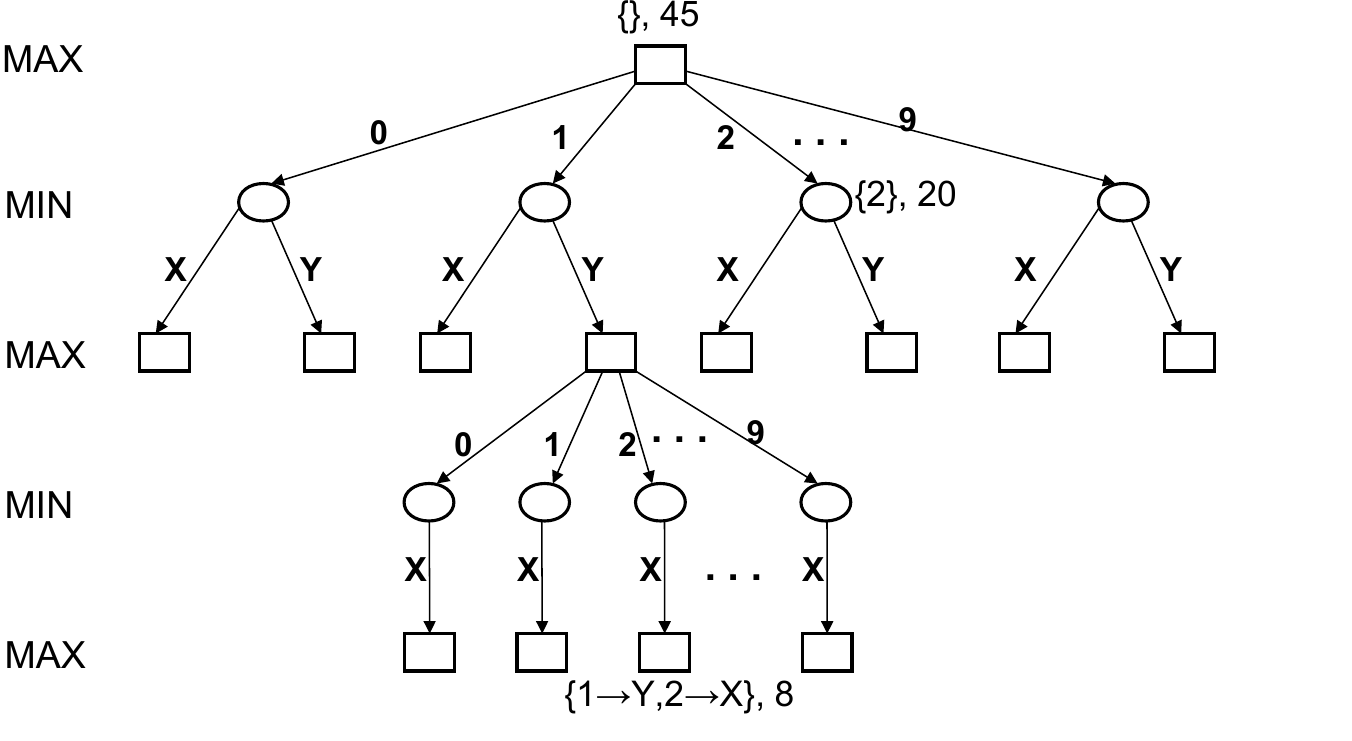} 
\caption{A partial game tree for $\E = (10-X)*Y$}
\end{center}
\label{fig:tree}
\end{figure}

\noindent
The game tree $\tree ((10-X)*Y)$ is shown partially in Figure~\ref{fig:tree}. The edges are labeled by the moves corresponding to the players.

The minimax value is computed by the method of backward induction applied to $\tree(\E)$. This procedure works by reasoning backwards from the end of the game and computing the optimal move for the players at each position. At a leaf node, the expression is a constant, and the value of the node is this constant. Working up, each MIN node has as its value, the minimum of the values of its valid children and each MAX node, the maximum of the values of its valid children\footnote{An internal node is deemed valid if it has a valid child. A leaf node is valid if it evaluates to a finite value.}. The value of the root is $\minimax (\E)$.
\smallskip

\section{Alpha-beta Pruning and Node Ordering}
\label{sec:alpha}

To determine $\minimax(\E)$, it's not necessary to evaluate every node in the tree. Suppose alpha is the current maximum (over the children evaluated so far) for a MAX node and beta the current minimum for a MIN node. For a MAX node, if its alpha value is at least the beta value of its parent, then there is no reason to explore the node further as the final value of its parent will be smaller than alpha. Each pruning of a subtree of a MAX node in this manner is referred to as a \emph{beta cutoff}. Similarly, for a MIN node, if its beta value is at most the alpha value of its parent then the remaining subtrees of this node can be pruned as the final value of its parent will be larger than beta. These prunings are \emph{alpha cutoffs}.
\smallskip

A subtree at height $d$  that is pruned by an alpha cutoff is rooted at a MAX node and prunes $T(d/2)$ nodes. A subtree at height $d$  that is pruned by a beta cutoff is rooted at a MIN node and prunes $(T(\lceil d/2 \rceil)-1)/10$ nodes. 
\smallskip

For example, suppose $\E = (10 - X)*Y$. Then $\minimax (\E) = 45$ and a leaf node that achieves this value is $\{5 \rightarrow X, 9 \rightarrow Y\}$. To compute $\minimax (\E)$, we start at the root node and evaluate the MIN nodes $\{0\}, \{1\}, \dots, \{5\}$ in succession which return the values $0, 10, 20, 30, 40$ and $45$ respectively. At this point, the alpha value at the root is $\max (0, 10, 20, 30, 40, 45) = 45$. When node $\{6\}$ is explored, the MIN node computes the value of the MAX node $\{6 \rightarrow X\}$ which returns $36$ as its minimax value. Thus the beta value of $\{6\}$ is $36$ which is smaller than $45$, the alpha value of its parent. As a result, the node $\{6 \rightarrow Y\}$ is not evaluated. Similarly, the nodes $\{7 \rightarrow Y\}, \{8 \rightarrow Y\}$ and $\{9 \rightarrow Y\}$ are not evaluated leading to 4 alpha cutoffs. 
\smallskip

The pseudocode for computing the minimax value of $\E$ with alpha-beta pruning is given by Algorithm~\ref{alg:1}. The function alphabeta() takes as its parameters, the current node, its height, the current value (alpha or beta) of the parent node, the digit passed (valid for a MIN node) and the current player. Apart from the minimax value, the algorithm also returns the number of nodes pruned by alpha and beta cutoffs, which are computed using the recursion formula (\ref{eqn:1}), as well as the entire principal variation. The function is called with the command \newline
\texttt{
alpha\_prunes = beta\_prunes = 0; principal\_var = ``'' \newline
alphabeta (root, $2n$, $\infty$, $-1$, MAX)
}
\smallskip

\begin{algorithm}
\caption{Minimax value of $\E$ with alpha-beta pruning}
\label{alg:1}
\begin{algorithmic}[0]
\Function{alphabeta}{node, height, parent\_$\alpha \beta$, digit, player}
	\If{height = 0} \Comment{leaf node}
		\State principal\_var = ``''
		\State \textbf{return} the value of node
	\EndIf
	\If{player = MAX}		\Comment{process MAX node}
		\State maxstr = ``''		\Comment{principal variation from this node}
		\State maxval = $-\infty$	\Comment{current $\alpha$}
		\For{each $i$ from 0 to 9}	\Comment{evaluate each child in this loop}
			\State value = alphabeta (node, height $ - 1$, maxval, $i$, MIN)
			\If{value $>$ maxval and value $\ne +\infty$}	\Comment{update $\alpha$}
				\State maxval = value
				\State maxstr = principal\_var + `i'
			\EndIf
			\If{maxval $\ge$ parent\_$\alpha \beta$}	\Comment{beta prune}
				\State beta\_prunes = beta\_prunes + $(9 - i)*(T(\mathrm{height}/2)-1)/10$
				\State break
			\EndIf
		\EndFor
		\State principal\_var = maxstr
		\State \textbf{return} maxval
	
	\Else				
		\State minstr = ``''		
		\State minval = $+\infty$	\Comment{current $\beta$}
		\State $j = (\mathrm{height}+1)/2$	\Comment{number of children left to explore}
		\For{each uninstantiated variable $v$ in node}
			\State j = j - 1
			\State grandchild = node (digit $\rightarrow v$)	\Comment{replace $v$ by digit in node}
			\State value = alphabeta (grandchild, $\mathrm{height} - 1$, minval, $-1$, MAX)
			\If{value $<$ minval and value $\ne -\infty$}	\Comment{update $\beta$}
				\State minval = value
				\State minstr = principal\_var + `v'
			\EndIf
			\If{minval $\le$ parent\_$\alpha \beta$}	\Comment{alpha prune}
				\State alpha\_prunes = alpha\_prunes + $j*T((\mathrm{height}-1)/2)$
				\State break
			\EndIf
		\EndFor
		\State principal\_var = minstr
		\State \textbf{return} minval
	\EndIf
			
\EndFunction
\end{algorithmic}
\end{algorithm}

\subsection{A Heuristic for Node Ordering}
\label{subsec:ordering}

The effectiveness of alpha-beta pruning depends on the order in which each node's children are explored. For example, with $\E = (10 - X)*Y$, suppose we evaluate a MAX node by choosing the digits in sequence $\{5, 4, 6, 3,$ $7, 2, 8, 1, 9, 0\}$, and evaluate a MIN node by setting the variable sequence as $(Y, X)$ if the digit passed to it is less than 5 and $(X, Y)$ otherwise. Then, to calculate $\minimax (\E)$, the node $\{5\}$ is evaluated first and returns $45$. Subsequently, for each of the MIN nodes, the order in which its children are explored ensures that there is an alpha cutoff. 
\smallskip

Let $\tilde{v} (x)$ be an estimate for the value $v(x)$ of node $x$. We propose a heuristic for determining the order in which the digits are to be chosen at a MAX node. The ordering is static in the sense that it is determined by $\E$ and is the same for all nodes being evaluated. We estimate the values of the MAX nodes 2 levels below the root node. These estimates are backed up, by taking the minima, to estimate the values of their parents. If these estimates are placed in decreasing order, as $\tilde{v}(\{i_0\}) \ge \tilde{v}(\{i_1\}) \ge \dots \ge \tilde{v}(\{i_9\})$ then the children of a MAX node are evaluated in sequence $i_0, i_1, \dots, i_9$.
\smallskip

For a MAX node $x = \{i \rightarrow X\}$, our estimate for $v(x)$ is simply the maximum of $\E$ over some random instantiations of the variables of $\E$ while fixing $X$ at $i$. More precisely, to estimate $v (\{i \rightarrow X\})$, we fix $X$ at $i$ and randomly instantiate the other variables in $\E$ with digits and compute $\val(\E)$. We do this a fixed number of times and take the maximum of the resulting values.
\smallskip

The performance of minimax algorithm is further enhanced by noting that  several nodes in $\tree (\E)$ correspond to the same game position and thus have to to be evaluated only once. An example are nodes $\{ 2 \rightarrow X, 1 \rightarrow Y \}$ and $\{1 \rightarrow Y, 2 \rightarrow X \}$ in $\tree (X*(Y - Z))$. We exploit this fact by storing, for each MAX node $x$ that is fully evaluated (i.e. none of its children are pruned by beta cutoffs), its value and the principal variation starting at $x$. On subsequent visits to nodes that correspond to the same game position, this value is retrieved instead of being recomputed.

\section{Implementation Details}
\label{sec:Imp}
We first convert $\E$ to a postfix form using Dijkstra's shunting yard algorithm \cite{Dijk}. During evaluation, the variables are substituted with values, and $\val(\E)$ is computed using the reverse polish evaluation \cite{BWW} algorithm.
\smallskip

Table~\ref{tab:1} compares the number of nodes pruned by alpha-beta and alpha-beta with node ordering and also shows the ordering of digits at each MAX node as determined by the heuristic. For the alpha-beta method, the number of nodes pruned is the sum of the number of nodes pruned by alpha and beta cutoffs. For alpha-beta with node ordering, the number of nodes pruned is the sum of the number of nodes pruned by alpha and beta cutoffs and the transposition tables. For expressions with five or six variables, we have observed a ten-fold speedup in the performance of the second method over the first.
\smallskip

We also attempted ordering the MIN nodes as well as using different orderings for MAX nodes at different heights using the same heuristic but any gains in the number of nodes pruned was offset by the computational time in determining the order. Other promising approaches such as Negascout \cite{Rei} and the MTD-\textit{f} \cite{PSPD} algorithm have not been attempted yet.

\begin{table}

\setlength{\extrarowheight}{4pt}
\setlength{\tabcolsep}{6pt} 
\begin{tabular}{|p{2.5cm}|p{2.5cm}|p{2cm}|p{2cm}|p{1.8cm}|}
\hline
\centering{\E} & $\minimax(\E)$ and principal variation & \small{$\alpha$-$\beta$: No. of nodes pruned} & \multicolumn{2}{|p{3.8cm}|}{$\alpha$-$\beta$ with node ordering and transposition tables}\\
\cline{4-5}
& & & Digit order & \centering{\small{No. of nodes pruned}} \tabularnewline
\cline{1-5} \cline{1-5}

\centering{$\frac{x}{y} + \frac{2*y}{z} - \frac{z}{x}$} & $16/3,\: 3 \rightarrow x,\:$ $3 \rightarrow z, \: 9 \rightarrow y$ & \centering{$4515$} & $1, 4, 3, 5, 6,$ $7, 2, 8, 9, 0$ & \centering{$7526$} \tabularnewline \cline{1-5}

$w - \frac{y*z}{3} + 3*x$ & $21, \: 6 \rightarrow w, 7 \rightarrow y,$ $5 \rightarrow x, \: 0 \rightarrow z$ & \centering{$302271$} & $9, 8, 7, 6, 5, $ $4, 3, 2, 1, 0$ & \centering{$464162$} \tabularnewline \cline{1-5}

$v + w + x - y - z$ & $12, \: 7 \rightarrow y, 8 \rightarrow v,$ $7 \rightarrow w, \: 4 \rightarrow x,$ $0 \rightarrow z$ &\centering{$ \approx 2.18*10^8$} & $7, 3, 4, 5, 6,$ $1, 2, 8, 9, 0$ & $\approx 2.53*10^8$ \\ \cline{1-5}

\centering{$\frac{a+b}{c} + \frac{d+e}{f}$} & $7.6,\: 3 \rightarrow a, 3 \rightarrow c,$ $5 \rightarrow f, \: 9 \rightarrow b,$ $9 \rightarrow d,$ $9 \rightarrow e$&\centering{$\approx 1.33*10^{10}$} & $5, 4, 3, 2, 8,$ $1, 9, 0, 6, 7$ & $\approx\! 1.55*10^{10}$ \\ \cline{1-5}

\end{tabular}
\caption{Comparison of alpha-beta and alpha-beta with node ordering}
\label{tab:1}
\end{table}

\section{e-Valuate on Graphs}	
\label{sec:G}
In this section, we describe two-player games on graphs where each player's goal is to maximize his share of a set $B$ of resources (which could be vertices, edges or faces -- for planar graphs). In each game, players MAX and MIN color the elements of a (different) set $A$ in the following manner. At each turn,  MAX proposes a color from a palette of available colors that MIN assigns to an uncolored element of $A$ of his choice. When all the elements of $A$ are colored, a criterion is used to determine which elements of $B$ belong to MAX and which to MIN.
\smallskip

The games we describe will be of the form $\game (G, A, B, [k], f)$ where
\begin{itemize}
\item $G$ is the graph on which the game is played. $G$ has vertex set $V$ and edge set $E$. When $G$ is planar, we will denote its set of faces by $R$.
\item $A$ is the set of elements to be colored.
\item $B$ is the set of elements to be shared.
\item $[k]$ is the set of colors and
\item $f$ prescribes the rule for determining whether a given $x \in B$ belongs to MAX or MIN. There is one such $f$ for each $x$ which we denote by $f_x$. $x$ belongs to MAX (resp. MIN) if $f_x$ evaluates to $1$ (resp. $0$). $f_x$ is defined in the following manner. We first associate a nonempty set $S(x) \subseteq A$ with $x$. $f_x$ evaluates to $0$ if all elements of $S(x)$ are assigned the same color and $1$ otherwise.
\end{itemize}

\noindent
At the conclusion of the game, MAX's share is given by evaluating the function $\F := \sum\limits_{x \in B} f_x$ on the input colors. 
\smallskip

If we interpret $\F$ as a mapping $\F : [k]^{|A|}\mapsto \N$ on the variable set $A$, then the game we have described is equivalent to playing e-Valuate on $\F$; MAX proposes values (colors) that MIN assigns to variables (elements of $A$) of his choice. Thus, under optimal play from both players, MAX's share is $\minimax(\F, [k])$. The winner can then be determined by comparing MAX's share with $\minimax(\F, [k])$.

\smallskip

If $|S(x)| > 1$ for each $x \in B$, then the game has the following monotonic property.
\begin{itemize}
\item If $k = 1$, then MAX has to choose the same color at each turn and ends up not acquiring any element of $B$. Thus $\minimax(\F,[1]) = 0$.
\item If $k = |A|$, then MAX chooses a different color at each turn and ends up acquiring all elements of $B$. Thus $\minimax(\F,[|A|]) = |B|$.
\item $\minimax(\F,[i]) \le \minimax(\F,[i+1])$ for $i = 1, 2, \dots, |A| - 1$.
\end{itemize}
This leads to the chain of inequalities:
$$0 = \minimax(\F, [1]) \le \minimax (\F, [2]) \le \dots \le \minimax(\F, [|A|]) = |B|$$
Thus an interesting problem is determining $\min \{k: \minimax(\F, [k]) = |B|\}$ for different classes of graphs.
\smallskip

The game is generalized by applying weights to each element of $B$. In this case, the goal of each player is to maximize the sum of the weights of his share. The corresponding function $\F$ is then a weighted sum of $f$s i.e. if $w(x)$ is the weight of $x \in B$, then $\F = \sum_{x \in B} w(x) f_x$.
\smallskip

There are six possible choices for the pair $(A, B)$: $(V, E)$, $(E, V)$, $(E, R)$, $(V, R)$, $(R, E)$ and $(R, V)$. We now describe the games corresponding to each pair. 
\smallskip

\subsection*{Sharing of Edges}	\label{subsec:GE}
This is $\game (G, V, E, [k], f)$, i.e. the set of edges are shared by coloring the vertices. For an edge $e = (u, v)$, let $S(e) = \{u, v\}$. Then $e$ belongs to MAX if $u$ and $v$ are assigned different colors and to MIN otherwise. 
\smallskip

\subsection*{Sharing of Vertices}	\label{subsec:GV}
This is $\game (G, E, V, [k], f)$, i.e. the set of vertices are shared by coloring the edges. For $v \in V$, $S(v) := \{e \in E: e \mbox{ is incident to }v\}$. As $S(v) \ne \emptyset$ we require that all vertices have positive degree. 
Then $v$ belongs to MAX if there are two edges incident to $v$ colored differently and to MIN otherwise. 

\smallskip

\subsection*{Sharing of Faces}	\label{subsec:GF}
Here we assume that $G$ is planar.
This is then $\game (G, E, R, [k], f)$, i.e. the set of faces are shared by coloring the edges. With a given face $r$, we associate the set of edges on its boundary, i.e. $S(r) = \{ e \in E: e \mbox{ lies on } r\mbox{'s boundary}\}$. This implies that $r$ belongs to MAX if there are two edges on its boundary colored differently and to MIN otherwise.
\smallskip

\noindent
We observe that
\begin{itemize}
\item Coloring the vertices instead of the edges leads to $\game (G, V, R, [k], f)$. The criterion for determining if a face belongs to MAX is similar; two vertices on its boundary should be assigned different colors. It's not clear which version is more favorable to MAX. 
\item A natural choice for the weights of the faces are their areas (we exclude the outer face by assigning it a weight of zero). In this weighted version, each player tries to maximize his share of the area enclosed by the graph.
\end{itemize}

Finally there are two other games involving coloring of faces to share vertices or edges: $\game (G, R, V, [k], f)$ and $\game (G, R, E, [k], f)$. If $x$ is a vertex or an edge, we associate with it the set of faces it lies on, i.e. $S(x) = \{r \in R: x \mbox{ lies on } r\mbox{'s boundary}\}$.
\medskip

\noindent
Here's an example demonstrating these games. 
\vspace*{0.5cm}

\begin{figure}[htbp]
\begin{center}
\includegraphics[height=3in, width=3in]{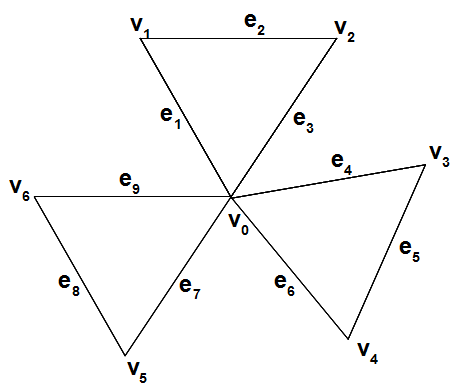} 
\caption{The Friendship Graph $F_3$}
\label{fig:G}
\end{center}
\end{figure}

\noindent
{\bf Example 1}: We illustrate the games we have described on the Friendship graph $F_3$ (see Figure~\ref{fig:G}). We assume $k = 2$ and optimal play from both players. For convenience, we assume that the colors are red and blue. We describe the edge-sharing game in detail. The other games use similar reasoning and we outline the steps. 
\smallskip

\noindent
{\game$\mathbf{(F_3, V, E, [2], f)}$}: We show that MAX's share is $s = 4$ edges. MAX ensures that $s \ge 4$ by proposing three red and four blue colors (in any order). If vertex $v_0$ is colored red, then it is adjacent to 4 vertices colored blue implying $s \ge 4$. If $v_0$ is colored blue, then it is adjacent to 3 vertices colored red. Moreover, at least one of the edges $(v_1, v_2)$, $(v_3, v_4)$, $(v_5, v_6)$ has its ends colored differently resulting in $s \ge 4$. 

On the other hand MIN ensures that $s \le 4$ by assigning red colors to the vertices $v_1, v_2, v_3, v_0, v_4, v_5, v_6$ (in that order) and blue colors in the reverse order (i.e. $v_6, v_5, v_4, v_0, v_3, v_2, v_1$). After all the vertices are colored, at least one of the 3 triangles will be monochromatic. The other two will each have utmost two edges whose ends are colored differently.

The function $\F$, such that $\minimax(\F, [2]) = 4$, is given by
\begin{eqnarray*}
\F & = & f(v_0, v_1) + f(v_1, v_2) + f(v_2, v_0) + f(v_0, v_3) + f(v_3, v_4) + \\ 
& & f(v_4, v_0) + f(v_0, v_5) + f(v_5, v_6) + f(v_6, v_0)
\end{eqnarray*}
where $f(x, y) = 1$ if $x \ne y$ and $0$ otherwise. 
\vspace*{0.3cm}

\noindent
{\game$\mathbf{(F_3, E, V, [2], f)}$}: Here MAX's share is $s = 2$ vertices. MAX chooses one red and eight blue colors (in any order) ensuring $s \ge 2$. MIN assigns red colors to the edges $e_1, \dots, e_9$ in that order and blue colors in the reverse order. This ensures $s \le 2$.
\vspace*{0.3cm}

\noindent
{\game$\mathbf{(F_3, E, R, [2], f)}$}: For this game, MAX's share is $2$ faces. The strategy for choosing and assigning colors to edges is same as in \game$(F_3, E, V, [2], f)$. Using this strategy, MIN ensures that two of the triangular faces are monochromatic. Thus MAX acquires one triangular face and the outer face.
\vspace*{0.3cm}

\noindent
{\game$\mathbf{(F_3, V, R, [2], f)}$}: For this game, MAX's share is $3$ faces. The strategy for choosing and assigning colors to vertices is same as in $\game(F_3, V, E, [2], f)$.
\vspace*{0.3cm}

\noindent
{\game$\mathbf{(F_3, R, V, [2], f)}$}: Let the three triangular faces and the outer face be labeled $T_1, T_2, T_3$ and $T_4$ respectively. MAX's share is $5$ vertices. MAX chooses two red and two blue colors. If $T_4$ is colored red (say), then the five vertices in the two triangular faces colored blue are acquired by MAX. MIN applies red colors to the faces $T_1, T_4, T_2, T_3$ (in that order) and blue colors in the reverse order ensuring that $T_4$ and one of the triangular faces share the same color.
\vspace*{0.3cm}

\noindent
{\game$\mathbf{(F_3, R, E, [2], f)}$}: MAX's share is $6$ edges. The strategy of both players is same as in \game${(F_3, R, V, [2], f)}$.
\qed

\medskip

Finally we note that other variants of the game can be derived by changing the rule $f$ that determines whether elements of $B$ belong to MAX or MIN. One variant when $k = 2$, is to assign an $x \in B$ to MAX if an odd number of elements in $S(x)$ are colored $1$. If MAX proposes color $1$ at every turn, then he will acquire all $x \in B$ such that $|S(x)|$ is odd (for example, in $\game (G, E, V, [2], f)$, MAX would acquire all vertices of odd degree). Thus care must be taken to ensure that $|S(x)|$ is even for each $x \in B$ by modifying $G$ if necessary. In the vertex sharing game this can be achieved by joining pairs of odd-degree vertices by edges. The resulting graph may not be simple but the game is still playable.

The interesting feature of this variant is that under random assignment of colors to $A$, an element in $B$ is equally likely to belong to MAX or MIN (since $S(x)$ has an equal number of subsets of odd and even size). Then by linearity of expectation, the expected share of MAX (and MIN) is $|B|/2$.

\section{A Tiling Game}
\label{sec:T}
In the standard two-player sequential game, the game tree is used to determine the player (if any) with the winning strategy. At the leaf node, the game has concluded and the value assigned to the node depends on the winner; usually $\pm 1$ or $0$ depending on whether the winner is MAX or MIN or it's a draw. These values are then backed up to the root to identify the player with the winning strategy.
\smallskip

On the other hand, in e-Valuate, the value at the leaf node is the evaluation of the function associated with the game. The winner is determined by comparing this value with that at the root of the tree (found by backing up values as before). Thus in e-Valuate, under optimal play, the game is \emph{always} a draw.
\smallskip

Thus once we have an activity that can be completed in different ways (such as coloring a graph, instantiating a function or vector, creating a Hamiltonian tour in a weighted graph etc.) and for each completion, we can associate a value, we have a game on; MAX and MIN complete the activity together but with contrasting goals of maximizing and minimizing the final value.
\smallskip

Here is a tiling game based on this idea.
\smallskip

\begin{figure}[htbp]
\begin{center}
\includegraphics[height=2.2in]{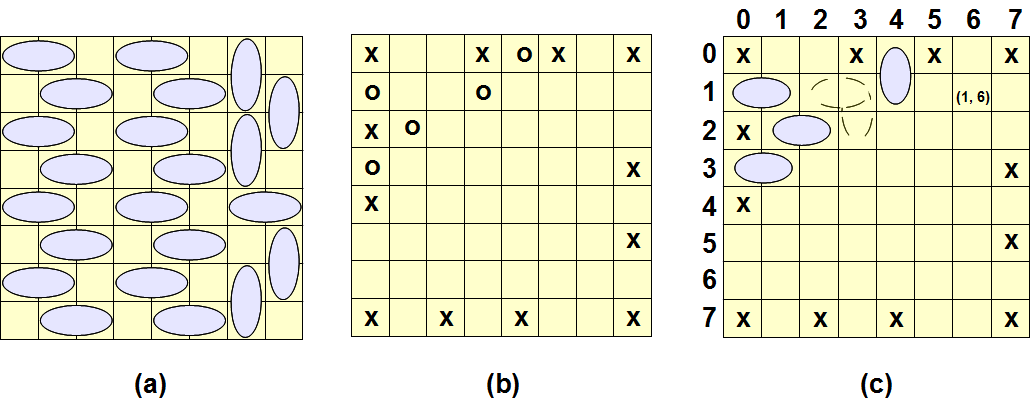} 
\caption{Tiling by dominoes. (a) A maximal partial tiling by $22$ dominoes, (b) A configuration with $4$ unoccupied squares along each border, (c) A partial tiling of (b)}
\label{fig:D}
\end{center}
\end{figure}

A partial tiling of an $8\times 8$ checkerboard is the placement of non-overlapping dominoes on the board such that each domino occupies exactly two squares (so each domino is either horizontal or vertical). The tiling is said to be \emph{maximal} if it cannot be extended any further, i.e. all unoccupied regions are $1 \times 1$ squares. Figure \ref{fig:D}(a) shows a maximal partial tiling with $22$ dominoes.
\smallskip

In the \emph{tiling game}, MAX and MIN take turns in placing dominoes on the board until the tiling is maximal. The value at the end of the game is the number of dominoes on the board. Thus, left to himself, MAX would like to cover every square of the board using $32$ dominoes. On the other hand, MIN would choose a maximal tiling that uses only $22$ dominoes as shown in Figure \ref{fig:D}(a). This is the minimum as proved in Theorem \ref{thm:D}. Thus under optimal play, the final number of dominoes $x$, satisfies $22 \le x \le 32$.
\smallskip

\begin{thm}
Any maximal partial tiling of an $8 \times 8$ checkerboard by dominoes requires at least $22$ dominoes.
\label{thm:D}
\end{thm}

\noindent
{\bf Proof}:
We label the squares of the board by their coordinates $(i, j): 0  < i, j < 8$ (see Figure~\ref{fig:D}(c)). A line refers to a row or column of the board. The maximality of the partial tiling implies the following easily verifiable facts:
\begin{enumerate}
\item Each line contains utmost $4$ unoccupied squares.
\item Each line intersects at least $3$ dominoes.
\end{enumerate}
For $j > 0$, we map an unoccupied square at $(i, j)$, to the domino on its left, i.e. occupying the square at $(i, j-1)$. This defines a one-to-one mapping from the set of unoccupied squares in columns $1$ to $7$ to the the set of dominoes not intersecting column $7$. If there are utmost $3$ unoccupied squares in column $0$, then these can be mapped to the $3$ dominoes intersecting column $7$. This completes the one-to-one mapping from the set of unoccupied squares to the set of dominoes. Thus if there are $x$ unoccupied squares, then number of dominoes $\ge x$ and together they occupy at least $3x$ squares implying $3x \le 64$. Since $x$ is even, this implies $x \le 20$ resulting in at least $22$ dominoes on the board.
\smallskip

A similar reasoning shows that if the rightmost column or the top or bottom row have utmost $3$ unoccupied squares, then there are at least $22$ dominoes.
\smallskip

We now consider the case where the lines along the border have $4$ unoccupied squares each. Our goal is to find an additional domino that is not mapped to any unoccupied square under the mapping described above. Then this, along with the three dominoes intersecting column $7$ can be mapped to the $4$ unoccupied squares on column $0$ and our proof will be complete.
\smallskip

A maximal partial tiling that has $4$ unoccupied squares along each of the border lines is possible only if the following conditions are satisfied.
\begin{itemize}
\item The $4$ corner squares are unoccupied.
\item If the square at $(0, 2)$ is unoccupied, then the square at $(2, 0)$ is occupied (since otherwise the dominoes occupying the squares at $(0, 1)$ and $(1, 0)$ would overlap) and vice versa.
\end{itemize}

These observations imply, apart from symmetry, a unique configuration of the unoccupied squares along the border lines as shown in Figure \ref{fig:D}(b); the unoccupied squares are marked by {\bf x} while some of the squares that must be occupied are marked by {\bf o}.
\smallskip

If we cover the squares marked by {\bf o} by dominoes, we get a partial tiling as shown in Figure \ref{fig:D}(c). If the domino at $(1, 3)$ is horizontal then $(1, 2)$ is occupied and the domino at $(1, 0)$ is left unmapped. Similarly if the domino at $(1, 3)$ is vertical then $(2, 3)$ is occupied and the domino at $(2,1)$ is left unmapped. In either case we have uncovered a domino left unmapped and this completes our proof.
\qed

\section{Conclusion}
\label{sec:conc}
We have demonstrated the effectiveness of search algorithms for computing the minimax value of e-Valuate. Other heuristics for node values could yield more effective ordering of the nodes of the game tree and thus faster algorithms. We have also described some games on graphs based on e-Valuate.
\smallskip

One would also like to understand what expressions and associated domains constitute a fair game. A fair game is one where if MAX and MIN make their moves randomly, they have equal chances of winning. For example, if $\E$ has only $+$ and $*$ operators or is defined on one variable, then $\minimax(\E) = \max(\E)$ and MIN can never lose. 
\smallskip

On the games on graphs, some unresolved questions are:
\begin{itemize}
\item Computation of the minimax value: Can the graph structure be exploited to compute the value more efficiently than using the game tree approach for arithmetic expressions?
\item Estimates for the minimax value in terms of the graph parameters and the number of colors.
\item Estimate of the number of colors so that MAX and MIN have approximately the same share under optimal play.
\end{itemize}

\medskip

\noindent
{\bf Acknowledgment.} The authors thank Rohit Sabharwal and Aman Arora for helpful discussions.

\end{document}